\documentclass{article}%
\usepackage{amsmath}
\usepackage{amsfonts}
\usepackage{amssymb}
\usepackage{graphicx}%
\newtheorem{theorem}{Theorem}

\newtheorem{remark}[theorem]{Remark}

\newenvironment{proof}[1][Proof]{\noindent\textbf{#1.} }{\ \rule{0.5em}{0.5em}}
\begin{document}

\title{The relation between the counting function $N\left(  \lambda\right)  $ and the
heat kernel $K\left(  t\right)  $}
\author{BY WU-SHENG DAI and MI XIE \thanks{We are very indebted to Dr G. Zeitrauman
for his encouragement. This work is supported in part by NSF of China, under
Project No.10605013 and No.10675088.}}
\date{}
\maketitle

\begin{abstract}
For a given spectrum $\left\{  \lambda_{n}\right\}  $ of the Laplace operator
on a Riemannian manifold, in this paper, we present a relation between the
counting function $N\left(  \lambda\right)  $, the\ number of eigenvalues
(with multiplicity) smaller than $\lambda$, and the heat kernel $K\left(
t\right)  $, defined by $K\left(  t\right)  =\sum_{n}e^{-\lambda_{n}t}$.
Moreover, we also give an asymptotic formula for $N\left(  \lambda\right)  $
and discuss when $\lambda\rightarrow\infty$ in what cases $N\left(
\lambda\right)  =K\left(  1/\lambda\right)  $.

\end{abstract}

The relation between the spectrum of the Laplace operator on a Riemannian
manifold and the geometry of this Riemannian manifold is an important subject
\cite{Berard,Berger,Milnor,GWW}, and the problem of spectral asymptotics is
one of the central problems in the theory of partial differential operators
\cite{Ivrii}. For a given spectrum $\left\{  \lambda_{n}\right\}  $ of the
Laplace operator on a Riemannian manifold, one can in principle obtain the
counting function $N\left(  \lambda\right)  $, defined to be%
\begin{equation}
N\left(  \lambda\right)  =\text{the\ number of eigenvalues (with multiplicity)
of the Laplace operator smaller than }\lambda\text{,}%
\end{equation}
and the heat kernel, defined to be%
\begin{equation}
K\left(  t\right)  =\sum_{n}e^{-\lambda_{n}t}.
\end{equation}
One of the main problems is to seek the asymptotic expansions of $N\left(
\lambda\right)  $ and $K\left(  t\right)  $. Usually, it is relatively easy to
obtain the asymptotic expansion of the heat kernel $K\left(  t\right)  $.
Nevertheless, it is difficult to calculate the asymptotic expansion of the
counting function $N\left(  \lambda\right)  $ \cite{Berger}. The
Hardy-Littlewood-Karamata Tauberian theorem gives the first term of the
asymptotic expansion of $N\left(  \lambda\right)  $ \cite{Kac}, but does not
provide any information beyond the first-order term. In this paper, we point
out a relation between $N\left(  \lambda\right)  $ and $K\left(  t\right)  $.

\begin{theorem}%
\begin{equation}
K\left(  t\right)  =t\int_{0}^{\infty}N\left(  \lambda\right)  e^{-\lambda
t}d\lambda. \label{035}%
\end{equation}

\end{theorem}

\begin{proof}
The generalized Abel partial summation formula reads%
\begin{equation}
\sum_{u_{1}<\lambda_{n}\leq u_{2}}b\left(  n\right)  f\left(  \lambda
_{n}\right)  =B\left(  u_{2}\right)  f\left(  u_{2}\right)  -B\left(
u_{1}\right)  f\left(  u_{1}\right)  -\int_{u_{1}}^{u_{2}}B\left(  u\right)
f^{\prime}\left(  u\right)  du,\label{040}%
\end{equation}
where $\lambda_{i}\in\mathbb{R}$, $\lambda_{1}\leq\lambda_{2}\leq\cdots
\leq\lambda_{n}\leq\cdots$, and $\lim_{n\rightarrow\infty}\lambda_{n}=\infty$.
$f\left(  u\right)  $ is a continuously differentiable function on $\left[
u_{1},u_{2}\right]  $ $\left(  0\leq u_{1}<u_{2}\text{, }\lambda_{1}\leq
u_{2}\right)  $, $b\left(  n\right)  $ $\left(  n=1,2,3,\cdots\right)  $ are
arbitrary complex numbers, and $B\left(  u\right)  =\sum_{\lambda_{n}\leq
u}b\left(  n\right)  $. We apply the generalized Abel partial summation
formula, Eq. (\ref{040}), with $f\left(  u\right)  =e^{-u\left(
s-s_{0}\right)  }$ and $b\left(  n\right)  =a_{n}e^{-\lambda_{n}s_{0}}$, where
$s$, $s_{0}\in\mathbb{C}$. Then%
\begin{equation}
A\left(  u_{2},s\right)  -A\left(  u_{1},s\right)  =A\left(  u_{2}%
,s_{0}\right)  e^{-u_{2}\left(  s-s_{0}\right)  }-A\left(  u_{1},s_{0}\right)
e^{-u_{1}\left(  s-s_{0}\right)  }+\left(  s-s_{0}\right)  \int_{u_{1}}%
^{u_{2}}A\left(  u,s_{0}\right)  e^{-u\left(  s-s_{0}\right)  }du,\label{050}%
\end{equation}
where%
\begin{equation}
A\left(  u,s\right)  =\sum_{\lambda_{n}\leq u}a_{n}e^{-\lambda_{n}%
s}.\label{055}%
\end{equation}

Setting $a_{n}=1$ in Eq. (\ref{055}), we find%
\[
A\left(  \lambda,0\right)  =\sum_{\lambda_{n}\leq\lambda}1=N\left(
\lambda\right)  ,
\]
the counting function, and%
\[
A\left(  \infty,t\right)  =\sum_{n}e^{-\lambda_{n}t}=K\left(  t\right)  ,
\]
the heat kernel. By Eq. (\ref{055}), we also have $A\left(  0,t\right)  =0$.
Then, by Eq. (\ref{050}), we have%
\begin{equation}
K\left(  t\right)  =A\left(  \infty,t\right)  -A\left(  0,t\right)  =t\int
_{0}^{\infty}N\left(  \lambda\right)  e^{-\lambda t}d\lambda. \label{060}%
\end{equation}
This is just Eq. (\ref{035}).
\end{proof}

Furthermore we can also obtain the following theorem.

\begin{theorem}%
\begin{equation}
N\left(  \lambda\right)  =\frac{1}{2\pi i}\int_{c-i\infty}^{c+i\infty}K\left(
t\right)  \frac{e^{\lambda t}}{t}dt,\text{ \ \ }c>\lim_{n\rightarrow\infty
}\frac{\ln n}{\lambda_{n}}. \label{065}%
\end{equation}

\end{theorem}

\begin{proof}
By the Perron formula, we have%
\begin{equation}
\sum_{\mu_{n}<x}a_{n}=\frac{1}{2\pi i}\int_{c-i\infty}^{c+i\infty}f\left(
t\right)  \frac{x^{t}}{t}dt,\label{080}%
\end{equation}
where%
\begin{equation}
f\left(  s\right)  =\sum_{n=1}^{\infty}\frac{a_{n}}{\mu_{n}^{s}},\label{070}%
\end{equation}
and $c$ is a constant which is greater than the abscissa of absolute
convergence of the Dirichlet series $f\left(  s\right)  $. Setting
\[
a_{n}=1\text{ \ and }\mu_{n}=e^{\lambda_{n}}%
\]
in Eq. (\ref{070}), we obtain the heat kernel,%
\[
f\left(  t\right)  =\sum_{n=1}^{\infty}e^{-\lambda_{n}t}=K\left(  t\right)  .
\]
The abscissa of absolute convergence of $f\left(  t\right)  $ equals its
abscissa of convergence, equaling $\overline{\lim}_{n\rightarrow\infty}\ln
n/\lambda_{n}=\lim_{n\rightarrow\infty}\ln n/\lambda_{n}$. Thus, by Eq.
(\ref{080}), we have%
\[
N\left(  \lambda\right)  =\sum_{\lambda_{n}<\lambda}1=\frac{1}{2\pi i}%
\int_{c-i\infty}^{c+i\infty}K\left(  t\right)  \frac{e^{\lambda t}}{t}dt,
\]
and $c>\lim_{n\rightarrow\infty}\ln n/\lambda_{n}$. This proves the theorem.
\end{proof}

The above two theorems give the relation between the counting function
$N\left(  \lambda\right)  $ and the heat kernel $K\left(  t\right)  $.

One of the reasons why the counting function $N\left(  \lambda\right)
=\sum_{\lambda_{n}<\lambda}1$ is very difficult to calculate is that one often
encounters some unsolved problems in number theory when calculating $N\left(
\lambda\right)  $. For example, when calculating the counting function for the
spectrum of the Laplace operator on a tori, one encounters the Gauss circle
problem in number theory \cite{Berger}. In the following we will give an
asymptotic formula for $N\left(  \lambda\right)  $.

\begin{theorem}%
\begin{equation}
N\left(  \lambda\right)  =\sum_{n}\frac{1}{e^{\beta\left(  \lambda_{n}%
-\lambda\right)  }+1},\text{ \ \ }\left(  \beta\rightarrow\infty\right)  .
\label{100}%
\end{equation}

\end{theorem}

\begin{proof}
Observing that%
\[
\lim_{\beta\rightarrow\infty}\frac{1}{e^{\beta\left(  \lambda_{n}%
-\lambda\right)  }+1}=\left\{
\begin{array}
[c]{cc}%
1, & \text{ when }\lambda_{n}<\lambda,\\
0, & \text{ when }\lambda_{n}>\lambda,
\end{array}
\right.
\]
we have%
\[
\lim_{\beta\rightarrow\infty}\sum_{n}\frac{1}{e^{\beta\left(  \lambda
_{n}-\lambda\right)  }+1}=\sum_{\lambda_{n}<\lambda}1=N\left(  \lambda\right)
.
\]

\end{proof}

\begin{remark}
The asymptotic formula for $N\left(  \lambda\right)  $ given by Eq.
(\ref{100}) converts a partial sum ($\sum_{\lambda_{n}<\lambda}$) into a sum
over all possible values ($\sum_{\lambda_{n}<\infty}$). This will make the
calculation somewhat easy.
\end{remark}

In some cases the counting function approximately equals the heat kernel.

\begin{theorem}
Let $\rho\left(  \lambda\right)  $ be the number of eigenstates per unit
interval (the density of eigenstates). In the limit $\lambda\rightarrow\infty$
or $t\rightarrow0$,%
\begin{equation}
N\left(  \lambda\right)  =K\left(  \frac{1}{\lambda}\right)  \text{\ or
\ }N\left(  \frac{1}{t}\right)  =K\left(  t\right)  , \label{090}%
\end{equation}
when $\rho\left(  \lambda\right)  $ is a constant.
\end{theorem}

\begin{proof}
In the limit $\lambda\rightarrow\infty$ or $t\rightarrow0$, the summations can
be converted into integrals:%
\begin{align*}
N\left(  \lambda\right)   &  =\sum_{\lambda_{n}\leq\lambda}1=\int_{0}%
^{\lambda}\rho\left(  \lambda^{\prime}\right)  d\lambda^{\prime},\\
K\left(  t\right)   &  =\sum_{n}e^{-\lambda_{n}t}=\int_{0}^{\infty}\rho\left(
\lambda^{\prime}\right)  e^{-\lambda^{\prime}t}d\lambda^{\prime}.
\end{align*}
If $\rho\left(  \lambda\right)  =C$, where $C$ is a constant, then%
\begin{align*}
N\left(  \lambda\right)   &  =C\lambda,\\
K\left(  t\right)   &  =\frac{C}{t}.
\end{align*}
This proves the theorem.
\end{proof}

This is just the case that Weyl \cite{Weyl}, Pleijel \cite{Pleijel}, and Kac
\cite{Kac} discussed.

\vskip 1cm

{\footnotesize SCHOOL OF SCIENCE, TIANJIN UNIVERSITY, TIANJIN, P. R. CHINA }

{\footnotesize LIUHUI CENTER FOR APPLIED MATHEMATICS, NANKAI UNIVERSITY \&
TIANJIN UNIVERSITY, TIANJIN, P. R. CHINA }

{\footnotesize E-\textit{mail address}: daiwusheng@tju.edu.cn\newline }

{\footnotesize SCHOOL OF SCIENCE, TIANJIN UNIVERSITY, TIANJIN, P. R. CHINA }

{\footnotesize LIUHUI CENTER FOR APPLIED MATHEMATICS, NANKAI UNIVERSITY \&
TIANJIN UNIVERSITY, TIANJIN, P. R. CHINA }

{\footnotesize E-\textit{mail address}: xiemi@tju.edu.cn }

\end{document}